\documentclass[11pt]{article}

\usepackage[a4paper,margin=1in]{geometry}
\usepackage[T1]{fontenc}
\usepackage{lmodern}
\usepackage{microtype}
\usepackage{amsmath,amssymb,amsthm,mathtools}
\usepackage{enumitem}
\usepackage{hyperref}
 
 \allowdisplaybreaks
\makeatletter
\def\keywords{\xdef\@thefnmark{}\@footnotetext}
\makeatother
\allowdisplaybreaks
\numberwithin{equation}{section}
\newtheorem{theorem}{Theorem}[section]
\newtheorem{proposition}[theorem]{Proposition}
\newtheorem{lemma}[theorem]{Lemma}
\newtheorem{corollary}[theorem]{Corollary}
\theoremstyle{remark}
\newtheorem{remark}[theorem]{Remark}

\newcommand{\R}{\mathbb R}
\newcommand{\N}{\mathbb N}
\newcommand{\Q}{\mathbb Q}
\newcommand{\E}{\mathbb E}
\newcommand{\Pmu}{\mathbb P_{X_0}}
\newcommand{\Emu}{\mathbb E_{X_0}}
\newcommand{\1}{\mathbf 1}
\newcommand{\cM}{\mathcal M}
\newcommand{\cH}{\mathcal H}
\newcommand{\cL}{\mathcal L}
\newcommand{\wtM}{\widetilde{\mathcal M}}
\newcommand{\wtH}{\widetilde{\mathcal H}}
\newcommand{\dd}{\,\mathrm d}

\DeclareMathOperator{\sgn}{sgn}

\author{
Ziyi Chen$^{*}$  \quad   Jieliang Hong$^\dagger$  
}

\title{On the differentiability of local times of $(1+\beta)$-stable super-Brownian motion}

\date{{\small  {\it  
$^{*}$$^{\dagger}$Department of Mathematics, Southern University of Science and Technology, \\Shenzhen, China\\
\quad \\
    $^{*}$E-mail:  {\tt 12431004@mail.sustech.edu.cn} \\
$^\dagger$E-mail:  {\tt hongjl@sustech.edu.cn} 
  }
  }
  }

\begin{document}

 \maketitle
\keywords{This work was supported by the National Natural Science Foundation of China (Grant number: 12571150).}
\keywords{{\bf AMS 2020 subject classification:} 60J55; 60J68; 60G57}
\keywords{{\bf Key words and phrases:} stable super-Brownian motion; local time; spatial derivative; c\'adl\'ag modification; Tanaka formula.}
\begin{abstract}

We consider the local time $L(t,x)$ of $(1+\beta)$-stable super-Brownian motion with $0<\beta<1$.  In dimension one, when $X_0$ is atomless, we prove that $L(t,\cdot)$ is differentiable at every fixed deterministic point almost surely and is differentiable Lebesgue-a.e.  Nevertheless, almost surely it is not differentiable at every spatial point simultaneously.
\end{abstract}

\section{Introduction}

Let $M_F=M_F(\R^d)$ denote the space of finite measures on $(\R^d,\mathcal B(\R^d))$, equipped with the topology of weak convergence.  Let $X=(X_t)_{t\ge0}$ be a super-Brownian motion with $(1+\beta)$-stable branching mechanism, $0<\beta<1$, started from a deterministic measure $X_0\in M_F$.  The associated occupation measure is
\begin{equation*}
 Y_t(\cdot):=\int_0^t X_s(\cdot)\dd s.
\end{equation*}
Whenever $Y_t$ is absolutely continuous, we write
\[
 Y_t(\dd x)=L(t,x)\dd x
\]
and call $L$ the local time.  Fleischmann~\cite{Fle88} proved that the local time exists when $d<2+2/\beta$.  Mytnik and Perkins~\cite{MP03} showed that for $0<\beta<1$ it has a jointly continuous version in dimension $d=1$, whereas in dimensions $2\le d<2+2/\beta$ it is locally unbounded on every open set carrying positive occupation mass.  By contrast, in the finite-variance case $\beta=1$, Sugitani~\cite{Sug89} proved spatial differentiability in dimension one under suitable assumptions on the initial measure.

We study the spatial differentiability of $L(t,\cdot)$ for $0<\beta<1$.  The martingale term in the Tanaka formula has a natural spatial c\`adl\`ag version whose jumps are determined by the jumps of the superprocess.  This yields the one-sided derivatives of the local time and its set of nondifferentiability points.

Let $P_t$ be the semigroup of standard Brownian motion and set, for $\lambda>0$,
\begin{equation*}
 G_\lambda(u)=\int_0^\infty e^{-\lambda s}p_s(u)\dd s
 =\frac1{\sqrt{2\lambda}}e^{-\sqrt{2\lambda}|u|},
 \qquad
 g_\lambda(u)=-\sgn(u)e^{-\sqrt{2\lambda}|u|},
\end{equation*}
where $\sgn(0)=0$.  Write
\[
 G_\lambda^x(y)=G_\lambda(x-y),\qquad
 g_\lambda^x(y)=g_\lambda(x-y).
\]
The Tanaka formula of Mytnik and Xiang~\cite{MX04,Xiang05} gives, in dimension one,
\begin{equation*}
 L(t,x)=X_0(G_\lambda^x)-X_t(G_\lambda^x)
 +\lambda\int_0^tX_s(G_\lambda^x)\dd s+M_t(G_\lambda^x).
\end{equation*}
Define
\begin{equation*}
 Z^\lambda(t,x):=L(t,x)-X_0(G_\lambda^x)
\end{equation*}
and the deterministic-index field
\begin{equation*}
 H^\lambda(t,x):=-X_t(g_\lambda^x)
 +\lambda\int_0^tX_s(g_\lambda^x)\dd s+M_t(g_\lambda^x).
\end{equation*}
For each fixed deterministic $x$, $H^\lambda(t,x)$ is the $L^1$ limit of the difference quotients of $Z^\lambda(t,\cdot)$.

\begin{lemma}
\label{lem:intro-initial-convolution}
For every $X_0\in M_F(\R)$ and $\lambda>0$, the map $x\mapsto X_0(G_\lambda^x)$ is globally Lipschitz.  Moreover, for every $x\in\R$,
\begin{equation}
 D_+X_0(G_\lambda^x)=X_0(g_\lambda^x)-X_0(\{x\}),
 \qquad
 D_-X_0(G_\lambda^x)=X_0(g_\lambda^x)+X_0(\{x\}).
 \label{eq:intro-initial-one-sided}
\end{equation}
In particular, if $X_0$ is atomless, then $x\mapsto X_0(G_\lambda^x)$ is continuously differentiable with derivative $X_0(g_\lambda^x)$.
\end{lemma}

\begin{proof}
Since $|G_\lambda(u)-G_\lambda(v)|\le |u-v|$,
\[
 |X_0(G_\lambda^{x+h})-X_0(G_\lambda^x)|\le X_0(1)|h|.
\]
Separating the atom at $x$ and applying dominated convergence on $\R\setminus\{x\}$ gives \eqref{eq:intro-initial-one-sided}.  If $X_0$ is atomless, the two derivatives coincide; continuity of $x\mapsto X_0(g_\lambda^x)$ follows by dominated convergence.
\end{proof}

Every jump of $X$ has the form $\Delta X_s=r_s\delta_{z_s}$.  For $y\in\R$, write informally
\[
 R_t(y)=\sum_{0<s\leq t} r_s\1_{\{z_s=y\}}.
\]
The definition by jump-size truncation is given in Section~\ref{sec:canonical-versions}; there $R_t(y)<\infty$ is proved simultaneously for all $y$.

\begin{theorem}
\label{thm:intro-main}
Let $d=1$, $0<\beta<1$, $X_0\in M_F(\R)$ be deterministic, and fix $t,\lambda>0$.
\begin{enumerate}[label=\textup{(\roman*)},leftmargin=2.4em]
\item The deterministic-index field $x\mapsto H^\lambda(t,x)$ admits a canonical spatial c\`adl\`ag modification $\cH_t^\lambda$.  On a single probability-one event,
\begin{equation*}
 \Delta_y\cH_t^\lambda=-2R_t(y),\qquad y\in\R.
\end{equation*}
Moreover, every open set $U$ satisfying $Y_t(U)>0$ contains a point $y$ with $R_t(y)>0$.  Hence every fixed modification of $H^\lambda(t,\cdot)$ is discontinuous on such a set and has infinite H\"older seminorm of every positive order there.

\item With probability one, $Z^\lambda(t,\cdot)$ is locally absolutely continuous and
\begin{equation*}
 D_xZ^\lambda(t,x)=\cH_t^\lambda(x)
 \quad\text{for Lebesgue-a.e. }x\in\R.
\end{equation*}
For every fixed deterministic $x_0\in\R$,
\begin{equation*}
 \Pmu\bigl(D_xZ^\lambda(t,x_0)=H^\lambda(t,x_0)\bigr)=1.
\end{equation*}
Both one-sided derivatives exist simultaneously at every spatial point and satisfy
\begin{equation*}
 D_+Z^\lambda(t,y)-D_-Z^\lambda(t,y)=-2R_t(y).
\end{equation*}
Thus
\begin{equation*}
 \{y:Z^\lambda(t,\cdot)\text{ is not differentiable at }y\}
 =\{y:R_t(y)>0\}.
\end{equation*}

\item The local time $L(t,\cdot)$ is locally absolutely continuous.  On one probability-one event, both one-sided derivatives exist for every $y\in\R$ and
\begin{equation*}
 D_+L(t,y)-D_-L(t,y)
 =-2\bigl(R_t(y)+X_0(\{y\})\bigr).
\end{equation*}
Consequently,
\begin{equation*}
 \{y:L(t,\cdot)\text{ is not differentiable at }y\}
 =\{y:R_t(y)>0\}\cup\{y:X_0(\{y\})>0\}.
\end{equation*}
If $X_0(\{x_0\})=0$ for a fixed deterministic $x_0$, then
\begin{equation*}
 \Pmu\left(
 D_xL(t,x_0)=H^\lambda(t,x_0)+X_0(g_\lambda^{x_0})
 \right)=1.
\end{equation*}
\end{enumerate}
If $X_0(1)>0$, then almost surely neither $Z^\lambda(t,\cdot)$ nor $L(t,\cdot)$ is differentiable at every spatial point simultaneously.
\end{theorem}

For each fixed deterministic $y$, $R_t(y)=0$ almost surely, although $R_t$ may be positive at random spatial points.

Section~\ref{sec:preliminaries} gives the fixed-index estimates.  Sections~\ref{sec:canonical-versions}--\ref{sec:occupied-open-sets} construct the c\`adl\`ag derivative field and identify its jumps.  Section~\ref{sec:differentiability} proves the differentiability results.  The cumulative c\`adl\`ag construction is given in Appendix~\ref{sec:appendix-cumulative}.

\medskip
\noindent\textbf{Convention on constants.} A constant of the form $C(a,b,\ldots)$ may depend on the displayed parameters. Constants whose values are unimportant and may change from line to line are denoted by $C,c,c_d,\ldots$.

\section{Moment calculations and fixed-index estimates}
\label{sec:preliminaries}

Throughout the remainder of the paper, fix
\begin{equation}
 d=1,\qquad \beta\in(0,1),\qquad t>0,\qquad \lambda>0,
 \label{eq:parameters}
\end{equation}
and let $X$ start from a deterministic $X_0\in M_F(\R)$.  We use the standard $M_F(\R)$-valued c\`adl\`ag realization with the usual augmented natural filtration; see Dawson~\cite[Section~6]{Daw93}.  Write $\mu(\phi)=\int\phi\dd\mu$. 
The branching mechanism has the L\'evy--Khintchine representation
\begin{equation*}
 \Psi(u)=u^{1+\beta}
 =c_\beta\int_0^\infty(e^{-ru}-1+ru)r^{-2-\beta}\dd r,
 \qquad
 c_\beta:=\frac{\beta(1+\beta)}{\Gamma(1-\beta)}.
\end{equation*}
Let $P_s$ and $p_s$ be the semigroup and transition density of standard one-dimensional Brownian motion.  For $0\leq u\leq t$, define
\begin{equation*}
 Y_u(B):=\int_0^uX_s(B)\dd s,
 \qquad B\in\mathcal B(\R).
\end{equation*}
The mean-measure identities are
\begin{equation}
 \Emu[X_s(\phi)]=X_0(P_s\phi),
 \qquad
 \Emu[Y_t(\phi)]=\int_0^tX_0(P_s\phi)\dd s
 \label{eq:mean-measures}
\end{equation}
for nonnegative measurable $\phi$.  In particular,
\begin{equation*}
 \Emu[Y_t(\R)]=tX_0(1)<\infty.
\end{equation*}

By~\cite[Theorem~1.3(a)]{MP03}, there is a probability-one event on which
\begin{equation}
 Y_t(\dd y)=L(t,y)\dd y
 \label{eq:occupation-density}
\end{equation}
and $y\mapsto L(t,y)$ is continuous.  On the same event, simultaneously for every compact $K\subset\R$ and every $y\in\R$,
\begin{equation}
 \sup_{z\in K}L(t,z)<\infty,
 \qquad
 Y_t(\{y\})=0.
 \label{eq:occupation-consequences}
\end{equation}
Moreover, by~\cite[Theorem~1.1(a)]{MP03}, $X_t$ has a continuous density almost surely and is therefore atomless.

Every jump of $X$ is of the form $\Delta X_s=r\delta_z$; see \cite[Lemma~1.6(a)--(b)]{FMW10}.  Let
\[
 N=\sum_{(s,z,r):\Delta X_s=r\delta_z}\delta_{(s,z,r)}.
\]
Its compensator is
\begin{equation}
 \widehat N=\widehat N(\dd s,\dd z,\dd r)
 =c_\beta\dd s\,X_{s-}(\dd z)r^{-2-\beta}\dd r.
 \label{eq:jump-compensator}
\end{equation}
Write $\widetilde N=N-\widehat N$.  Since $X_{s-}=X_s$ for Lebesgue-a.e. $s$, either may be used in $ds$-integrals.  For $\phi\in C_b^2(\R)$, the martingale term is
\begin{equation*}
 M_t(\phi) :=\int_{(0,t]\times\R}\phi(z) \,M(ds,dz) :=\int_{(0,t]\times\R\times(0,\infty)} r\phi(z)\,\widetilde N(\dd s,\dd z,\dd r).
\end{equation*}
The stochastic integral is with respect to the martingale measure $M(ds,dz)$; see \cite[Section~6]{Daw93}.  For $s>0$, set
\[
(P_sX_0)(y):=\int_{\R}p_s(y-z)X_0(\dd z).
\]
For $p\geq 1$, let
 \begin{align*}
\cL_{loc}^p:=\cL_{loc}^p(\R^+ \times \R,   \, P_sX_0 (z) ds dz)
\end{align*}
denote the space of equivalence classes of measurable functions $\psi: \R^+ \times \R \to \R$ such that 
 \begin{align*}
 \int_0^T ds \int_{\R} |\psi(s,z)|^p (P_sX_0) (z) dz <\infty, \quad \forall T>0.
\end{align*}

For $p\in(1+\beta,2)$ and $\psi\in\cL_{loc}^p$, \cite[Lemma~1.7]{FMW10} implies that
\begin{align*}
 u\mapsto M_u(\psi):=\int_{(0,u]\times \R} \psi(s,z) M(ds, dz)  
\end{align*}
is well-defined.  The following estimate is from \cite[Lemma~3.1]{LM05}; see also \cite[Lemma~2.6]{FMW10}.

\begin{lemma} \label{stable-martingale-estimate}
For any $p\in (1+\beta, 2)$, $q\in (1, 1+\beta)$ and $T>0$, there exists a constant $C=C(\beta, p, q, T)>0$ such that for  any $\psi \in \cL_{loc}^p$,
\begin{align*}
\E_{X_0}\Big[\sup_{0\leq u\leq T} |M_u(\psi)|^q\Big] \leq & C \int_0^T ds \int_{\R} |\psi(s,y)|^q \,  (P_sX_0)(y) dy\\
&+C\Big[ \int_0^T ds \int_{\R} |\psi(s,y)|^p (P_sX_0) (y) dy\Big]^{q/p}.
\end{align*}
\end{lemma}

Recall the kernels
\begin{equation*}
 G_\lambda(u)=\frac1{\sqrt{2\lambda}}e^{-\sqrt{2\lambda}|u|},
 \qquad
 g_\lambda(u)=-\sgn(u)e^{-\sqrt{2\lambda}|u|},
\end{equation*}
with $\sgn(0)=0$, and write $G_\lambda^x(y)=G_\lambda(x-y)$ and $g_\lambda^x(y)=g_\lambda(x-y)$.  The Tanaka formula of Mytnik and Xiang~\cite[Theorem~1.1]{MX04} gives
\begin{equation}
 L(t,x)=X_0(G_\lambda^x)-X_t(G_\lambda^x)
        +\lambda Y_t(G_\lambda^x)+M_t(G_\lambda^x).
 \label{eq:tanaka-formula}
\end{equation}
Define
\begin{equation*}
 Z^\lambda(t,x):=L(t,x)-X_0(G_\lambda^x)
\end{equation*}
and the deterministic-index field
\begin{equation}
 H^\lambda(t,x)
 :=-X_t(g_\lambda^x)+\lambda Y_t(g_\lambda^x)+M_t(g_\lambda^x).
 \label{eq:def-H}
\end{equation}

\begin{proposition}
\label{prop:fixed-index-estimates}
For every $x\in\R$,
\begin{equation}
 \lim_{h\to0}\Emu\left|
 \frac{Z^\lambda(t,x+h)-Z^\lambda(t,x)}h-H^\lambda(t,x)
 \right|=0.
 \label{eq:L1-difference-quotient}
\end{equation}
Moreover, for every $K>0$ and $p\in(1+\beta,2)$ there is a constant $C=C(K,p,X_0,t,\lambda)>0$ such that
\begin{equation}
 \Emu|H^\lambda(t,x)-H^\lambda(t,z)|
 \leq C|x-z|^{1/p},
 \qquad x,z\in[-K,K].
 \label{eq:H-L1-modulus}
\end{equation}
In particular, $\Emu|H^\lambda(t,x)|<\infty$ for every $x\in\R$.
\end{proposition}

\begin{proof}
For fixed $x\in\R$ and $h\neq0$, set
\[
 \eta_h^x:=\frac{G_\lambda^{x+h}-G_\lambda^x}{h}-g_\lambda^x.
\]
Since $G_\lambda$ is $1$-Lipschitz and differentiable on $\R\setminus\{0\}$ with derivative $g_\lambda$,
\[
 \eta_h^x(y)\longrightarrow0\quad\text{for every }y\neq x,
 \qquad
 \|\eta_h^x\|_\infty\leq2.
\]
By \eqref{eq:tanaka-formula} and \eqref{eq:def-H},
\[
 \frac{Z^\lambda(t,x+h)-Z^\lambda(t,x)}{h}-H^\lambda(t,x)
 =-X_t(\eta_h^x)+\lambda Y_t(\eta_h^x)+M_t(\eta_h^x).
\]
The mean-measure identities and dominated convergence give
\[
 \Emu [X_t(|\eta_h^x|)]
 =\int_\R|\eta_h^x(y)|(P_tX_0)(y)\dd y\longrightarrow0
\]
and
\[
 \Emu [Y_t(|\eta_h^x|)]
 =\int_0^t\!\int_\R|\eta_h^x(y)|(P_sX_0)(y)\dd y\dd s
 \longrightarrow0.
\]
Fix $q\in(1,1+\beta)$ and $p\in(1+\beta,2)$.  Applying Lemma \ref{stable-martingale-estimate} to the time-independent integrand $\eta_h^x$, the same dominated-convergence argument yields
\[
 \Emu|M_t(\eta_h^x)|^q\longrightarrow0.
\]
Consequently,
\[
\begin{aligned}
 &\Emu\left|
 \frac{Z^\lambda(t,x+h)-Z^\lambda(t,x)}{h}-H^\lambda(t,x)
 \right|\\
 &\qquad\leq
 \Emu [X_t(|\eta_h^x|)]
 +\lambda\Emu [Y_t(|\eta_h^x|)]
 +\bigl(\Emu|M_t(\eta_h^x)|^q\bigr)^{1/q}
 \longrightarrow0,
\end{aligned}
\]
which proves \eqref{eq:L1-difference-quotient}.

For \eqref{eq:H-L1-modulus}, fix $x<z$ in $[-K,K]$ and choose $q\in(1,1+\beta)$.  On $\R\setminus[x,z]$, the two arguments of $g_\lambda$ have the same sign, so
\begin{equation}
 |g_\lambda^x(y)-g_\lambda^z(y)|
 \leq \sqrt{2\lambda}|x-z|,
 \qquad y\notin[x,z],
 \label{eq:g-outside-bound}
\end{equation}
whereas the global bound is $|g_\lambda^x-g_\lambda^z|\leq2$.  Also,
\begin{equation}
 \int_0^t\!\int_x^z(P_sX_0)(y)\dd y\dd s
 \leq C_tX_0(1)|x-z|.
 \label{eq:heat-content-bound}
\end{equation}
Apply Lemma \ref{stable-martingale-estimate} separately on $[x,z]$ and its complement.  With $\delta=z-x$, equations \eqref{eq:g-outside-bound} and
\eqref{eq:heat-content-bound} give
\begin{equation*}
 \Emu|M_t(g_\lambda^x)-M_t(g_\lambda^z)|^q
 \leq C(\delta^q+\delta+\delta^{q/p})
 \leq C_{K,p,q}|x-z|^{q/p}.
\end{equation*}
Therefore, by H\"older's inequality,
\begin{equation}
 \Emu|M_t(g_\lambda^x)-M_t(g_\lambda^z)|
 \leq C_{K,p}|x-z|^{1/p}.
 \label{eq:M-L1-increment}
\end{equation}
For the nonmartingale terms, the same decomposition and \eqref{eq:mean-measures} give, for example,
\[
 \Emu [X_t(|g_\lambda^x-g_\lambda^z|)]
 \leq \sqrt{2\lambda}\delta X_0(1)
      +2X_0(P_t\1_{[x,z]})
 \leq C_{t,\lambda,X_0}\delta,
\]
and integration in time yields the analogous occupation-measure estimate. Hence
\begin{align}
 \Emu|X_t(g_\lambda^x)-X_t(g_\lambda^z)|
 &\leq C_{K,t,\lambda,X_0}|x-z|,\notag\\
 \Emu|Y_t(g_\lambda^x)-Y_t(g_\lambda^z)|
 &\leq C_{K,t,\lambda,X_0}|x-z|.
 \label{eq:XY-L1-increments}
\end{align}
Combining \eqref{eq:def-H}, \eqref{eq:M-L1-increment}, and \eqref{eq:XY-L1-increments}, and using $|x-z|\leq2K$, proves \eqref{eq:H-L1-modulus}.  Finiteness of $\Emu|H^\lambda(t,x)|$ follows from the same estimate and Lemma \ref{stable-martingale-estimate}.
\end{proof}

\section{Spatial structure of the derivative field}
\label{sec:canonical-versions}

For $n\ge1$ and $y\in\R$, set
\begin{equation*}
 R_t^{(n)}(y):=
 \int_{(0,t]\times\{y\}\times(2^{-n},\infty)}
 r\,N(\dd s,\dd z,\dd r).
\end{equation*}
For each $n$ this is pathwise finite, and $R_t^{(n)}(y)$ is nondecreasing in $n$.  Define
\begin{equation}
 R_t(y):=\lim_{n\to\infty}R_t^{(n)}(y)\in[0,\infty].
 \label{eq:site-jump-mass}
\end{equation}
Equivalently,
\[
 R_t(y)=\int_{(0,t]\times\{y\}\times(0,\infty)}
 r\,N(\dd s,\dd z,\dd r)
\]
as an extended nonnegative integral.  Lemma~\ref{lem:cumulative-cadlag} proves simultaneous finiteness.

\begin{remark}

For fixed deterministic $y\in\R$, the compensation formula and $Y_t(\{y\})=0$ give
\[
 \Emu N((0,t]\times\{y\}\times(2^{-n},\infty))=0,\qquad n\ge1.
\]
Thus $R_t(y)=0$ almost surely for each deterministic $y$; nonzero values may occur at random sites.
\end{remark}

\subsection{Canonical c\`adl\`ag modifications and jump identities}

\begin{lemma}
\label{lem:cumulative-cadlag}
Define the deterministic-index field
\begin{equation*}
 \overline A_t(x):=
 \begin{cases}
  M_t(\1_{(0,x]}),&x\geq0,\\
  -M_t(\1_{(x,0]}),&x<0.
 \end{cases}
\end{equation*}
There is a spatial c\`adl\`ag process $A_t=(A_t(x))_{x\in\R}$ such that, for every fixed deterministic $x\in\R$,
\begin{equation}
 A_t(x)=\overline A_t(x)\qquad\Pmu\text{-a.s.}
 \label{eq:cumulative-modification}
\end{equation}
Moreover, there is a single event $\Omega_A$ with $\Pmu(\Omega_A)=1$ such that, on $\Omega_A$, simultaneously for every $y\in\R$,
\begin{equation}
 A_t(y)-A_t(y-)=R_t(y)<\infty.
 \label{eq:cumulative-fibre-jumps}
\end{equation}
\end{lemma}

\begin{proof}
See Appendix~\ref{sec:appendix-cumulative}.  Its Step~2 uses the independent counting estimate in Lemma~\ref{lem:counting-bounds}, proved in Section~\ref{sec:occupied-open-sets}; no conclusion of the present lemma is used in that estimate.
\end{proof}

For a spatial c\`adl\`ag function $f$, write $\Delta_yf:=f(y)-f(y-)$.  On the event $\Omega_A$ supplied by Lemma~\ref{lem:cumulative-cadlag},
\begin{equation}
 \Delta_yA_t=R_t(y),\qquad y\in\R,
 \label{eq:cumulative-jumps}
\end{equation}
simultaneously for every $y\in\R$.

Decompose
\begin{equation}
 g_\lambda=\sigma+\rho_\lambda,
 \qquad
 \sigma(u):=-\sgn(u),
 \qquad
 \rho_\lambda(u):=\sgn(u)(1-e^{-\sqrt{2\lambda}|u|}).
 \label{eq:kernel-decomposition}
\end{equation}
Then $\rho_\lambda$ is bounded and globally Lipschitz:
\begin{equation}
 |\rho_\lambda(u)-\rho_\lambda(v)|
 \leq \sqrt{2\lambda}|u-v|,
 \qquad u,v\in\R.
 \label{eq:rho-lipschitz}
\end{equation}

\begin{lemma}
\label{lem:lipschitz-part}
The deterministic-index field $x\mapsto M_t(\rho_\lambda^x)$ has a continuous modification on $\R$, where $\rho_\lambda^x(y)=\rho_\lambda(x-y)$. The modification is chosen jointly measurable in $(\omega, x)$.
\end{lemma}

\begin{proof}
Fix $q\in(1,1+\beta)$ and $p\in(1+\beta,2)$. Lemma \ref{stable-martingale-estimate} gives
\begin{align}
 \E|M_t(\rho_\lambda^x)-M_t(\rho_\lambda^z)|^q
 &\leq C\int_0^t\!\int_\R
   |\rho_\lambda(x-y)-\rho_\lambda(z-y)|^q
   (P_sX_0)(y)\dd y\dd s \notag\\
 &\quad+C\left(\int_0^t\!\int_\R
   |\rho_\lambda(x-y)-\rho_\lambda(z-y)|^p
   (P_sX_0)(y)\dd y\dd s\right)^{q/p},
 \label{eq:stable-moment-rho}
\end{align}
where $(P_sX_0)(y)=\int_\R p_s(u-y)X_0(\dd u)$. 
Applying \eqref{eq:rho-lipschitz} and $\int_\R(P_sX_0)(y)\dd y=X_0(1)$ to \eqref{eq:stable-moment-rho}, we get
\begin{equation*}
 \Emu|M_t(\rho_\lambda^x)-M_t(\rho_\lambda^z)|^q
 \leq C|x-z|^q.
\end{equation*}
Since $q>1$, by Kolmogorov's continuity theorem, the field $M_t(\rho_\lambda^x)$ admits a jointly measurable continuous modification.
\end{proof}

Fix the continuous version of $M_t(\rho_\lambda^x)$ from Lemma~\ref{lem:lipschitz-part}.
For $x\in\R$, put $\sigma^x(y):=\sigma(x-y)$.  If $x>0$, then
\begin{equation}
 \sigma^x-\sigma^0
 =-2\1_{(0,x)}-\1_{\{0\}}-\1_{\{x\}},
 \label{eq:sign-identity-positive}
\end{equation}
whereas, if $x<0$,
\begin{equation}
 \sigma^x-\sigma^0
 =2\1_{(x,0)}+\1_{\{x\}}+\1_{\{0\}}.
 \label{eq:sign-identity-negative}
\end{equation}

For every fixed deterministic $a\in\R$,
\begin{equation}
 M_t(\1_{\{a\}})=0\qquad\Pmu\text{-a.s.}
 \label{eq:deterministic-singleton-zero}
\end{equation}
By \eqref{eq:jump-compensator} and \eqref{eq:occupation-consequences},
\[
 \int_{(0,t]\times\{a\}\times(1,\infty)}r\,\widehat N(\dd s,\dd z,\dd r)
 =c_\beta Y_t(\{a\})\int_1^\infty r^{-1-\beta}\dd r=0 \qquad\Pmu\text{-a.s.}
\]
 By the compensation formula and nonnegativity of $N$, 
\begin{align*}
 \int_{(0,t]\times\{a\}\times(1,\infty)}r\, N(\dd s,\dd z,\dd r)=0 \qquad\Pmu\text{-a.s.}
\end{align*}
For the small-jump part, the martingale-measure isometry gives
\begin{align*}
 \Emu\left|\int_{(0,t]\times\{a\}\times(0,1]}
 r\,\widetilde N(\dd s,\dd z,\dd r)\right|^2
 &=\Emu\int_{(0,t]\times\{a\}\times(0,1]}
 r^2\,\widehat N(\dd s,\dd z,\dd r)\\
 &=c_\beta\Emu[Y_t(\{a\})]\int_0^1r^{-\beta}\dd r=0.
\end{align*}
Thus \eqref{eq:deterministic-singleton-zero} holds.  Therefore \eqref{eq:sign-identity-positive}--\eqref{eq:sign-identity-negative} and Lemma~\ref{lem:cumulative-cadlag} imply
\begin{equation}
 M_t(\sigma^x)=M_t(\sigma^0)-2A_t(x)
 \qquad\Pmu\text{-a.s. for every fixed deterministic }x.
 \label{eq:sign-martingale-identity}
\end{equation}

Define
\begin{equation*}
 \cM_t^\lambda(x)
 :=M_t(\rho_\lambda^x)+M_t(\sigma^0)-2A_t(x),
 \qquad x\in\R.
\end{equation*}
Then $\cM_t^\lambda$ is spatially c\`adl\`ag and, by
\eqref{eq:kernel-decomposition} and \eqref{eq:sign-martingale-identity},
\begin{equation}
 \cM_t^\lambda(x)=M_t(g_\lambda^x)
 \qquad\Pmu\text{-a.s. for every fixed deterministic }x.
 \label{eq:M-modification}
\end{equation}

The fields $A_t$ and $M_t(\rho_\lambda^\cdot)$ are jointly measurable, hence so is $\cM_t^\lambda$.  Any two c\`adl\`ag modifications are spatially indistinguishable.

\begin{lemma}
\label{lem:canonical-fields}
The field
\begin{equation}
 \cH_t^\lambda(x)
 :=-X_t(g_\lambda^x)+\lambda Y_t(g_\lambda^x)+\cM_t^\lambda(x)
 \label{eq:canonical-H}
\end{equation}
is a jointly measurable spatial c\`adl\`ag modification of $x\mapsto H^\lambda(t,x)$.  With probability one, simultaneously for every $y\in\R$,
\begin{equation}
 \Delta_y\cM_t^\lambda
 =\Delta_y\cH_t^\lambda
 =-2R_t(y).
 \label{eq:common-jumps}
\end{equation}
\end{lemma}

\begin{proof}
If $\mu$ is a finite atomless measure and $x_n\to x$, then $g_\lambda(x_n-y)\to g_\lambda(x-y)$ for $\mu$-a.e. $y$ and $|g_\lambda|\leq1$.  Dominated convergence therefore makes $x\mapsto\mu(g_\lambda^x)$ continuous.  Apply this pathwise to $X_t$ and $Y_t$, using their atomlessness from Section~\ref{sec:preliminaries}.  The two nonmartingale terms in \eqref{eq:canonical-H} are continuous and jointly measurable by measurable kernel integration, while \eqref{eq:M-modification} shows that $\cH_t^\lambda$ is a modification of $H^\lambda(t,\cdot)$.  Finally, Lemma~\ref{lem:lipschitz-part} and \eqref{eq:cumulative-jumps} give $\Delta_y\cM_t^\lambda=-2R_t(y)$ simultaneously for all $y$; adding the continuous nonmartingale terms gives the same identity for $\cH_t^\lambda$.
\end{proof}

\subsection{Discontinuity on occupied open sets}
\label{sec:occupied-open-sets}

\begin{lemma}
\label{lem:counting-bounds}
Fix $t>0$.  Let $(N_u)_{0\leq u\leq t}$ be a nonexplosive simple adapted counting process with $N_0=0$, and let $(\Lambda_u)_{0\leq u\leq t}$ be its continuous predictable compensator, with $\Lambda_0=0$ and $\Lambda_t<\infty$ almost surely.  For $a,v>0$,
\begin{equation}
 \mathbb P(N_t=0,\ \Lambda_t\geq a)\leq e^{-a},
 \label{eq:no-jump-bound}
\end{equation}
and, for a universal constant $C<\infty$,
\begin{equation}
 \mathbb P(N_t\geq2,\ \Lambda_t\leq v)\leq Cv^2.
 \label{eq:two-jump-bound}
\end{equation}
\end{lemma}

\begin{proof}
For each $\theta\in\R$,
\[
 \exp\{\theta N_u-(e^\theta-1)\Lambda_u\}
\]
is a nonnegative local martingale and hence, after localization and Fatou's lemma, a supermartingale. Hence, for $0\le u\le t$, 
$$\E[\exp\{\theta N_u-(e^\theta-1)\Lambda_u\}]\le1.$$
For $\theta<0$, its value on $\{N_t=0,\Lambda_t\geq a\}$ is at least $\exp\{(1-e^\theta)a\}$.  Markov's inequality followed by $\theta\downarrow-\infty$ proves \eqref{eq:no-jump-bound}.  For $\theta>0$,
\[
 \mathbb P(N_t\geq2,\Lambda_t\leq v)
 \leq\exp\{-2\theta+(e^\theta-1)v\}.
\]
If $0<v<1$, take $e^\theta=v^{-1}$; if $v\geq1$, use the trivial bound by one.  Thus \eqref{eq:two-jump-bound} holds, for instance with $C=e$.
\end{proof}

\begin{lemma}
\label{lem:occupied-jump}
For every fixed nonempty open interval $I\subset\R$,
\begin{equation*}
 \{Y_t(I)>0\}
 \subseteq\{\exists y\in I:R_t(y)>0\}
 \qquad\Pmu\text{-a.s.}
\end{equation*}
\end{lemma}

\begin{proof}
For $n\geq1$, a Borel set $B\subset\R$, and $0\leq u\leq t$, define
\begin{equation*}
 N_u^{(n)}(B):=N((0,u]\times B\times(1/n,\infty))
\end{equation*}
and its continuous compensator
\begin{equation*}
 \Lambda_u^{(n)}(B)
 :=c_\beta\int_0^uX_s(B)\dd s
   \int_{1/n}^\infty r^{-2-\beta}\dd r
 =\kappa_\beta n^{1+\beta}Y_u(B),
\end{equation*}
where $\kappa_\beta:=c_\beta/(1+\beta)$.  Each jump time contributes at most one mark to $N^{(n)}(I)$, and
\[
 \Emu[N_t^{(n)}(I)]=\Emu[\Lambda_t^{(n)}(I)]<\infty.
\]
Thus $N^{(n)}(I)$ is a simple, nonexplosive counting process, so Lemma~\ref{lem:counting-bounds} applies. By \eqref{eq:no-jump-bound}, for $m,n\geq1$,
\begin{equation}
 \Pmu\bigl(Y_t(I)>m^{-1},\ N_t^{(n)}(I)=0\bigr)
 \leq\exp\{-\kappa_\beta n^{1+\beta}/m\}.
 \label{eq:occupied-no-jump}
\end{equation}
Set
\[
 B_n:=(0,t]\times I\times(1/n,\infty),
 \qquad
 B:=(0,t]\times I\times(0,\infty).
\]
Since $B_n\uparrow B$, continuity from below of the pathwise measure $N(\omega;\cdot)$ gives
\[
 N(B)=\lim_{n\to\infty}N(B_n)
     =\lim_{n\to\infty}N_t^{(n)}(I).
\]
By \eqref{eq:site-jump-mass} and strict positivity of every jump size,
\begin{equation*}
 \{N(B)=0\}
 =\bigcap_{n\geq1}\{N_t^{(n)}(I)=0\}
 =\{R_t(y)=0\text{ for every }y\in I\}.
\end{equation*}
For $m\geq1$, let
\begin{equation*}
 A_m:=\{Y_t(I)>m^{-1}\}\cap\{N(B)=0\}.
\end{equation*}
For any $n\ge1$, $A_m$ is contained in the event on the left of \eqref{eq:occupied-no-jump}.  Letting $n\to\infty$ therefore gives $\Pmu(A_m)=0, ~m\ge1$.  Since
\[
 \{Y_t(I)>0,\ R_t(y)=0\text{ for every }y\in I\}
 =\bigcup_{m\geq1}A_m,
\]
the result follows.
\end{proof}

For an open $U\subset\R$ and $\gamma>0$, write
\begin{equation*}
 [f]_{\gamma;U}
 :=\sup_{\substack{x,z\in U\\x\neq z}}
 \frac{|f(x)-f(z)|}{|x-z|^\gamma}.
\end{equation*}

\begin{theorem}\label{thm:spatial-discontinuity} 
There is an event $\Omega_{\mathrm{sp}}$ with $\Pmu(\Omega_{\mathrm{sp}})=1$ such that, on this event, simultaneously for every open $U\subset\R$,
\begin{equation}
 Y_t(U)>0
 \quad\Longrightarrow\quad
 \exists y\in U:\quad
 \Delta_y\cM_t^\lambda
 =\Delta_y\cH_t^\lambda
 =-2R_t(y)\neq0.
 \label{eq:simultaneous-open-jump}
\end{equation}
On the same event, simultaneously for every open $U\subset\R$ and every
$\gamma>0$,
\begin{equation}
 Y_t(U)>0
 \quad\Longrightarrow\quad
 [\cM_t^\lambda]_{\gamma;U}
 =[\cH_t^\lambda]_{\gamma;U}=\infty.
 \label{eq:canonical-holder-failure}
\end{equation}

Moreover, let $\wtM_t^\lambda$ and $\wtH_t^\lambda$ be any fixed modifications of the two deterministic-index fields.  Then, almost surely, simultaneously for every open $U\subset\R$ with $Y_t(U)>0$, neither $\wtM_t^\lambda|_U$ nor $\wtH_t^\lambda|_U$ is continuous, and
\begin{equation}
 [\wtM_t^\lambda]_{\gamma;U}
 =[\wtH_t^\lambda]_{\gamma;U}=\infty
 \qquad\text{for every }\gamma>0.
 \label{eq:arbitrary-modification-failure}
\end{equation}
\end{theorem}

\begin{proof}
Let
\[
 \mathcal I_\Q:=\{(a,b):a,b\in\Q,\ a<b\}.
\]
This is a countable basis for the topology of $\R$. For $I\in\mathcal I_\Q$, define
\[
 E_I:=\{Y_t(I)=0\}\cup\{\exists y\in I:R_t(y)>0\}.
\]
Lemma~\ref{lem:occupied-jump} gives $\Pmu(E_I)=1$.  Intersecting over the
countable family $\mathcal I_\Q$, let $\Omega_{\mathrm{jump}}$ be the full-probability event from Lemma~\ref{lem:canonical-fields} on which both canonical paths are c\`adl\`ag and \eqref{eq:common-jumps} holds for every $y$, and set
\begin{equation*}
 \Omega_{\mathrm{sp}}
 :=\Omega_{\mathrm{jump}}\cap
   \bigcap_{I\in\mathcal I_\Q}E_I,
 \qquad \Pmu(\Omega_{\mathrm{sp}})=1.
\end{equation*}

Fix $\omega\in\Omega_{\mathrm{sp}}$ and an open $U\subset\R$ with $Y_t(U)>0$.  Since
\begin{equation*}
 U=\bigcup_{\substack{I\in\mathcal I_\Q\\I\subset U}}I,
\end{equation*}
there is $I\in\mathcal I_\Q$ with $I\subset U$ and $Y_t(I)>0$; otherwise
\[
 Y_t(U)\leq
 \sum_{\substack{I\in\mathcal I_\Q\\I\subset U}}Y_t(I)=0,
\]
a contradiction. Hence $R_t(y)>0$ for some $y\in I$, and \eqref{eq:common-jumps} proves \eqref{eq:simultaneous-open-jump}.  Since \(U\) is open, choose \(x_n\in U\) such that \(x_n<y\) and
\(x_n\uparrow y\).  By the existence of the left limit,
\[
 \bigl|\cH_t^\lambda(y)-\cH_t^\lambda(x_n)\bigr|
 \longrightarrow
 \bigl|\cH_t^\lambda(y)-\cH_t^\lambda(y-)\bigr|
 =
 \bigl|\Delta_y\cH_t^\lambda\bigr|>0,
\]
whereas $|y-x_n|^\gamma\to0$.  Thus $[\cH_t^\lambda]_{\gamma;U}=\infty$, and the same argument applies to $\cM_t^\lambda$.  This gives \eqref{eq:canonical-holder-failure}.

Let \(\wtM_t^\lambda\) and \(\wtH_t^\lambda\) be arbitrary fixed modifications.  For every fixed
\(q\in\Q\), the modification property gives
\[
 \wtM_t^\lambda(q)=\cM_t^\lambda(q),
 \qquad
 \wtH_t^\lambda(q)=\cH_t^\lambda(q)
 \qquad \Pmu\text{-a.s.}
\]
Since \(\Q\) is countable, the event
\begin{equation*}
 \Omega_{\Q}
 :=
 \bigcap_{q\in\Q}
 \left\{
 \wtM_t^\lambda(q)=\cM_t^\lambda(q),\
 \wtH_t^\lambda(q)=\cH_t^\lambda(q)
 \right\}
\end{equation*}
has probability one. Fix \(\omega\in\Omega_{\mathrm{sp}}\cap\Omega_{\Q}\), and let \(U\subset\R\) be open with
\(Y_t(U)>0\).  By \eqref{eq:simultaneous-open-jump}, there exists
\(y\in U\) such that
\begin{equation}
 \Delta_y\cM_t^\lambda = \Delta_y\cH_t^\lambda = -2R_t(y)\neq0.
\label{eq:open-jump}
\end{equation}
 Choose rational sequences $q_n^-\uparrow y, ~q_n^+\downarrow y, ~q_n^\pm\in\Q\cap U$. The c\`adl\`ag property gives
 \[
    \cH_t^\lambda(q_n^-) \to \cH_t^\lambda(y-), \qquad \cH_t^\lambda(q_n^+)\to\cH_t^\lambda(y).
 \]
The two limits are distinct by \eqref{eq:open-jump}, while on $\Omega_\Q$
\[
 \wtH_t^\lambda(q_n^-)=\cH_t^\lambda(q_n^-),
 \qquad
 \wtH_t^\lambda(q_n^+)=\cH_t^\lambda(q_n^+), \qquad n\ge1.
\]
Therefore the two rational sequences for $\wtH_t^\lambda$ have distinct limits, so $\wtH_t^\lambda$ is not continuous on $U$.  Moreover,
\[
 |q_n^+-q_n^-|\longrightarrow0,
 \qquad
 |\wtH_t^\lambda(q_n^+)-\wtH_t^\lambda(q_n^-)|
 \longrightarrow |\Delta_y\cH_t^\lambda|>0,
\]
and hence $[\wtH_t^\lambda]_{\gamma;U}=\infty$ for every $\gamma>0$.  The same argument applies to $\wtM_t^\lambda$.  This gives \eqref{eq:arbitrary-modification-failure}.
\end{proof}

\begin{corollary}
\label{cor:no-global-continuous-modification}
If $X_0(1)>0$, then every fixed modification of either $x\mapsto M_t(g_\lambda^x)$ or $x\mapsto H^\lambda(t,x)$ almost surely has no continuous sample path on $\R$.
\end{corollary}

\begin{proof}
Fix $\omega$ outside a null set so that $s\mapsto X_s(1)$ is right continuous at zero. There is a $\delta(\omega)>0$ such that $X_s(1)\ge X_0(1)/2>0$  for any $s\in(0,t\wedge\delta(\omega))$. Therefore
\[
Y_t(\R) = \int_0^t X_s(1)\,ds \ge\frac{X_0(1)}2(t\wedge\delta)>0.
\]
 Apply Theorem~\ref{thm:spatial-discontinuity} with $U=\R$.
\end{proof}

\section{Spatial differentiability of the local time}
\label{sec:differentiability}

\begin{proposition}
\label{prop:fixed-point-continuity}
For every fixed deterministic $x_0\in\R$,
\begin{equation*}
 \Pmu(\cH_t^\lambda\text{ is continuous at }x_0)=1,
 \qquad
 \Pmu(\cM_t^\lambda\text{ is continuous at }x_0)=1.
\end{equation*}
\end{proposition}

\begin{proof}
Choose $K>|x_0|+1$, fix $p\in(1+\beta,2)$, and set $x_n=x_0-n^{-1}$.  Then
\begin{equation}
 \cH_t^\lambda(x_n)\longrightarrow\cH_t^\lambda(x_0-)
 \qquad\Pmu\text{-a.s.}
 \label{eq:left-limit-sequence}
\end{equation}
Since $\cH_t^\lambda$ is a modification of $H^\lambda(t,\cdot)$, Proposition~\ref{prop:fixed-index-estimates} gives
\begin{equation}
 \Emu|\cH_t^\lambda(x_n)-\cH_t^\lambda(x_0)|
 \leq C_{K,p}n^{-1/p}.
 \label{eq:left-limit-L1}
\end{equation}
Fatou's lemma, \eqref{eq:left-limit-sequence}, and \eqref{eq:left-limit-L1} yield
\begin{align*}
 \Emu|\cH_t^\lambda(x_0-)-\cH_t^\lambda(x_0)|
 &\leq\liminf_{n\to\infty}
 \Emu|\cH_t^\lambda(x_n)-\cH_t^\lambda(x_0)|=0.
\end{align*}
Thus the left limit equals the value almost surely, while right-continuity is part of the c\`adl\`ag property.  The first assertion follows.  The common-jump identity \eqref{eq:common-jumps} then gives $\Delta_{x_0}\cM_t^\lambda=\Delta_{x_0}\cH_t^\lambda=0$ almost surely, which proves the second assertion.
\end{proof}

\begin{remark}
The assertion is for each fixed deterministic $x_0$ and does not imply simultaneous continuity in $x$.
\end{remark}

\begin{lemma}
\label{lem:zero-derivative-process}
Let $(W(x))_{x\in\R}$ have continuous sample paths and satisfy $\E|W(x)|<\infty$ for every $x$.  If
\begin{equation*}
 \lim_{h\to0}\frac1{|h|}\E|W(x+h)-W(x)|=0,
 \qquad x\in\R,
\end{equation*}
then, with probability one,
\begin{equation*}
 W(x)=W(0),\qquad x\in\R.
\end{equation*}
\end{lemma}

\begin{proof}
Set $f(x)=\E|W(x)-W(0)|$.  The reverse triangle inequality gives
\[
 \frac{|f(x+h)-f(x)|}{|h|}
 \leq\frac{\E|W(x+h)-W(x)|}{|h|}\longrightarrow0.
\]
Hence $f'(x)=0$ for every $x$, so $f$ is constant and $f(x)=f(0)=0$.  Thus $W(x)=W(0)$ almost surely for each fixed $x$.  Taking a countable intersection over $x\in\Q$ and using sample-path continuity gives the simultaneous conclusion for every $x\in\R$.
\end{proof}

\begin{proposition}
\label{prop:integral-representation}
There is a probability-one event on which, simultaneously for every
$x\in\R$,
\begin{equation}
 Z^\lambda(t,x)
 =Z^\lambda(t,0)+\int_0^x\cH_t^\lambda(z)\dd z.
 \label{eq:integral-representation}
\end{equation}
\end{proposition}

\begin{proof}
Use the continuous version of $Z^\lambda(t,\cdot)$ and define
\begin{equation*}
 W(x):=Z^\lambda(t,x)-Z^\lambda(t,0)
       -\int_0^x\cH_t^\lambda(z)\dd z.
\end{equation*}
A c\`adl\`ag function is locally bounded, so the integral is pathwise well-defined and continuous; hence $W$ has continuous paths.  The Tanaka formula, boundedness of $G_\lambda$, \eqref{eq:mean-measures}, and Lemma \ref{stable-martingale-estimate} show that $\Emu|Z^\lambda(t,x)|<\infty$.  If $I_x$ is the compact interval with endpoints $0$ and $x$, then Proposition~\ref{prop:fixed-index-estimates} and Tonelli's theorem give
\[
 \Emu\left|\int_0^x\cH_t^\lambda(z)\dd z\right|
 \leq\int_{I_x}\Emu|H^\lambda(t,z)|\dd z<\infty.
\]
Thus $\Emu|W(x)|<\infty$ for every $x$.

Fix $x$, choose $K>|x|+1$ and $p\in(1+\beta,2)$, and let $I(x,h)$ be the closed interval with endpoints $x$ and $x+h$.  Because $\cH_t^\lambda$ is a deterministic-index modification, for $0<|h|<1$,
\begin{align*}
 \frac1{|h|}\Emu|W(x+h)-W(x)|
 &\leq\Emu\left|
 \frac{Z^\lambda(t,x+h)-Z^\lambda(t,x)}h-H^\lambda(t,x)
 \right|\notag\\
 &\quad+\frac1{|h|}\int_{I(x,h)}
 \Emu|\cH_t^\lambda(z)-\cH_t^\lambda(x)|\dd z\notag\\
 &\leq o(1)+\frac{C_{K,p}}{|h|}\int_0^{|h|}u^{1/p}\dd u
 =o(1)+\frac{C_{K,p}}{1+1/p}|h|^{1/p}.
\end{align*}
The right-hand side tends to zero by Proposition~\ref{prop:fixed-index-estimates}.  Lemma \ref{lem:zero-derivative-process} gives $W(x)=W(0)=0$ simultaneously for all $x$, proving \eqref{eq:integral-representation}.
\end{proof}

For a real-valued function $f$, write
\[
 D_+f(y):=\lim_{h\downarrow0}\frac{f(y+h)-f(y)}h,
 \qquad
 D_-f(y):=\lim_{h\downarrow0}\frac{f(y)-f(y-h)}h,
\]
whenever the limits exist.

\begin{theorem}
\label{thm:spatial-differentiability}
Fix $t,\lambda>0$ as in \eqref{eq:parameters}.
\begin{enumerate}[label=\textup{(\roman*)},leftmargin=2.4em]
\item With probability one, $x\mapsto Z^\lambda(t,x)$ is locally absolutely continuous and
\begin{equation*}
 D_xZ^\lambda(t,x)=\cH_t^\lambda(x)
 \quad\text{for Lebesgue-a.e. }x\in\R.
\end{equation*}

\item For every fixed deterministic $x_0\in\R$,
\begin{equation*}
 \Pmu\bigl(D_xZ^\lambda(t,x_0)=H^\lambda(t,x_0)\bigr)=1.
\end{equation*}

\item On one probability-one event, simultaneously for every open
$U\subset\R$,
\begin{equation}
 Y_t(U)>0\quad\Longrightarrow\quad
 \exists y\in U:\quad
 \begin{cases}
  D_-Z^\lambda(t,y)=\cH_t^\lambda(y-),\\
  D_+Z^\lambda(t,y)=\cH_t^\lambda(y),\\
  D_+Z^\lambda(t,y)-D_-Z^\lambda(t,y)=-2R_t(y)\neq0.
 \end{cases}
 \label{eq:random-corner}
\end{equation}
In particular, the two-sided derivative of $Z^\lambda(t,\cdot)$ does not exist at the point $y$ furnished by \eqref{eq:random-corner}.
\end{enumerate}
\end{theorem}

\begin{proof}
The c\`adl\`ag path $\cH_t^\lambda$ is locally bounded and hence locally integrable.  Assertion~\textup{(i)} follows from \eqref{eq:integral-representation} and the fundamental theorem for Lebesgue integrals.

Fix a deterministic $x_0$.  By Proposition~\ref{prop:fixed-point-continuity}, $\cH_t^\lambda$ is continuous at $x_0$ almost surely.  Therefore \eqref{eq:integral-representation} gives
\[
 \frac{Z^\lambda(t,x_0+h)-Z^\lambda(t,x_0)}h
 =\frac1h\int_{x_0}^{x_0+h}\cH_t^\lambda(z)\dd z
 \longrightarrow\cH_t^\lambda(x_0)=H^\lambda(t,x_0)
\]
almost surely, proving~\textup{(ii)}.

Finally, work on the intersection of the probability-one events in Proposition~\ref{prop:integral-representation} and Theorem~\ref{thm:spatial-discontinuity}.  Let $U$ be open with $Y_t(U)>0$ and choose $y$ as in \eqref{eq:simultaneous-open-jump}.  For $h>0$,
\[
 \frac{Z^\lambda(t,y+h)-Z^\lambda(t,y)}h
 =\frac1h\int_y^{y+h}\cH_t^\lambda(z)\dd z
 \longrightarrow\cH_t^\lambda(y)
\]
by right-continuity.  Similarly,
\[
 \frac{Z^\lambda(t,y)-Z^\lambda(t,y-h)}h
 =\frac1h\int_{y-h}^{y}\cH_t^\lambda(z)\dd z
 \longrightarrow\cH_t^\lambda(y-).
\]
The difference of these limits is $-2R_t(y)\neq0$ by \eqref{eq:common-jumps}, proving~\textup{(iii)}.
\end{proof}

\begin{corollary}
\label{cor:local-time}
For arbitrary deterministic $X_0\in M_F(\R)$, the map $x\mapsto L(t,x)$ is locally absolutely continuous almost surely.  Moreover, there is a probability-one event on which, simultaneously for every $y\in\R$, both one-sided derivatives of $L(t,\cdot)$ exist and
\begin{equation}
 D_+L(t,y)-D_-L(t,y)
 =-2\bigl(R_t(y)+X_0(\{y\})\bigr),
 \qquad y\in\R.
 \label{eq:local-time-corner}
\end{equation}
Consequently, on the same event,
\begin{equation}
 \{y\in\R:L(t,\cdot)\text{ is not differentiable at }y\}
 =\{y:R_t(y)>0\}\cup\{y:X_0(\{y\})>0\}.
 \label{eq:local-time-exceptional-set}
\end{equation}
Moreover, for every fixed deterministic $x_0$ satisfying $X_0(\{x_0\})=0$,
\begin{equation}
 \Pmu\bigl(D_xL(t,x_0)
 =H^\lambda(t,x_0)+X_0(g_\lambda^{x_0})\bigr)=1.
 \label{eq:fixed-local-time-derivative}
\end{equation}
In particular, the last conclusion holds at every fixed deterministic point when $X_0$ is atomless.
\end{corollary}

\begin{proof}
By Lemma~\ref{lem:intro-initial-convolution}, the map $x\mapsto X_0(G_\lambda^x)$ is globally Lipschitz and hence locally absolutely continuous.  Together with Theorem~\ref{thm:spatial-differentiability}\textup{(i)} and $L(t,x)=Z^\lambda(t,x)+X_0(G_\lambda^x)$, this proves the local absolute continuity of $L(t,\cdot)$.

The one-sided derivative formulas in Lemma~\ref{lem:intro-initial-convolution} give, for every $y\in\R$,
\begin{equation*}
 D_+X_0(G_\lambda^y)
 =X_0(g_\lambda^y)-X_0(\{y\}),
 \qquad
 D_-X_0(G_\lambda^y)
 =X_0(g_\lambda^y)+X_0(\{y\}).
\end{equation*}
On the intersection of the probability-one events in
Proposition~\ref{prop:integral-representation} and
Lemma~\ref{lem:canonical-fields}, right-continuity and the existence of left
limits of $\cH_t^\lambda$ give, simultaneously for every $y\in\R$,
\[
 D_+Z^\lambda(t,y)=\cH_t^\lambda(y),
 \qquad
 D_-Z^\lambda(t,y)=\cH_t^\lambda(y-).
\]
Since $L(t,y)=Z^\lambda(t,y)+X_0(G_\lambda^y)$, it follows that
\[
 \begin{aligned}
 D_+L(t,y)
 &=\cH_t^\lambda(y)+X_0(g_\lambda^y)-X_0(\{y\}),\\
 D_-L(t,y)
 &=\cH_t^\lambda(y-)+X_0(g_\lambda^y)+X_0(\{y\}).
 \end{aligned}
\]
Subtracting these identities and using \eqref{eq:common-jumps} proves
\eqref{eq:local-time-corner}.  Since $R_t(y)$ and $X_0(\{y\})$ are
nonnegative, the two one-sided derivatives agree exactly when
\[
 R_t(y)=0
 \qquad\text{and}\qquad
 X_0(\{y\})=0.
\]
Because both one-sided derivatives exist at every $y$, this proves
\eqref{eq:local-time-exceptional-set}.

Finally, fix a deterministic $x_0$ satisfying $X_0(\{x_0\})=0$.
Proposition~\ref{prop:fixed-point-continuity} gives
\[
 \cH_t^\lambda(x_0-)=\cH_t^\lambda(x_0)
 \qquad\Pmu\text{-a.s.},
\]
while the modification property gives
\[
 \cH_t^\lambda(x_0)=H^\lambda(t,x_0)
 \qquad\Pmu\text{-a.s.}
\]
Substitution in the preceding one-sided derivative identities shows that
the two derivatives agree and proves
\eqref{eq:fixed-local-time-derivative}.
\end{proof}

\begin{corollary}
\label{cor:global-nondifferentiability}
If $X_0(1)>0$, then
\begin{align*}
 \Pmu\bigl(Z^\lambda(t,\cdot)\text{ is differentiable at every }x\in\R\bigr)
 &=0,\notag\\
 \Pmu\bigl(L(t,\cdot)\text{ is differentiable at every }x\in\R\bigr)
 &=0.
\end{align*}
\end{corollary}

\begin{proof}
The proof of Corollary~\ref{cor:no-global-continuous-modification} gives $Y_t(\R)>0$ almost surely.  Apply Theorem~\ref{thm:spatial-differentiability}\textup{(iii)} with $U=\R$, and then use \eqref{eq:local-time-exceptional-set} for the local time.
\end{proof}

Theorem~\ref{thm:intro-main} now follows from Lemma~\ref{lem:canonical-fields}, Theorem~\ref{thm:spatial-discontinuity}, Theorem~\ref{thm:spatial-differentiability}, Corollary~\ref{cor:local-time}, and Corollary~\ref{cor:global-nondifferentiability}.

\section*{Author contributions}
Research was conducted as a collaborative effort.

\section*{Funding}
This work was supported by the National Natural Science Foundation of China under Grant No.~12571150.

\section*{Data availability}
No datasets were generated or analyzed during the current study.

\section*{Declarations}
\noindent\textbf{Competing interests.} The authors declare no competing interests.

\appendix

\section{Proof of the cumulative c\`adl\`ag lemma}
\label{sec:appendix-cumulative}

For $\varepsilon>0$ and bounded Borel $\phi$, set
\begin{align*}
 M_t^\varepsilon(\phi)
 &:={}
 \int_{(0,t]\times\R\times(\varepsilon,\infty)}
 r\phi(z)\,N(\dd s,\dd z,\dd r)\notag\\
 &\quad-
 \int_{(0,t]\times\R\times(\varepsilon,\infty)}
 r\phi(z)\,\widehat N(\dd s,\dd z,\dd r).
\end{align*}
Both terms are finite almost surely, since
\begin{equation}
 \Emu N((0,t]\times\R\times(\varepsilon,\infty))
 =\frac{c_\beta}{1+\beta}\varepsilon^{-1-\beta}
   \Emu[Y_t(\R)]<\infty,
 \label{eq:truncated-count-integrability}
\end{equation}
so the first integral is a finite sum.  Moreover,
\begin{align}
 &\int_{(0,t]\times\R\times(\varepsilon,\infty)}
 r|\phi(z)|\,\widehat N(\dd s,\dd z,\dd r)\notag\\
 &\qquad\leq
 \frac{c_\beta}{\beta}\varepsilon^{-\beta}
 \|\phi\|_\infty Y_t(\R)<\infty
 \label{eq:truncated-compensator-integrability}
\end{align}
almost surely.

Set $\varepsilon_n=2^{-n}$ and
\begin{equation}
 A_t^{\varepsilon_n}(x):=
 \begin{cases}
  M_t^{\varepsilon_n}(\1_{(0,x]}),&x\geq0,\\
  -M_t^{\varepsilon_n}(\1_{(x,0]}),&x<0.
 \end{cases}
 \label{eq:truncated-cumulative}
\end{equation}
For each $n$, the raw-jump part is a finite step function and the compensator part is continuous; for $x\geq0$ it equals
\[
 \frac{c_\beta}{\beta}\varepsilon_n^{-\beta}
 \int_0^xL(t,z)\dd z.
\]
Thus every $A_t^{\varepsilon_n}$ is c\`adl\`ag and of finite variation.  Its dyadic increments are square-integrable.

We use the classical Rademacher--Menshov inequality: if $\xi_1,\ldots,\xi_N$ are orthogonal in $L^2$ and $S_k=\sum_{j=1}^k\xi_j$, then
\begin{equation}
 \E\max_{1\leq k\leq N}|S_k|^2
 \leq C\log^2(N+1)\sum_{j=1}^N\E[\xi_j^2].
 \label{eq:rademacher-menshov}
\end{equation}
See~\cite[Proposition~2.1]{Meaney07}.

For $n\geq1$, let
\begin{equation*}
 D_n(x):=A_t^{\varepsilon_{n+1}}(x)-A_t^{\varepsilon_n}(x).
\end{equation*}
We prove that, for every $m\in\N$,
\begin{equation}
 \sum_{n=1}^\infty\sup_{x\in[-m,m]}|D_n(x)|<\infty
 \qquad\Pmu\text{-a.s.}
 \label{eq:summable-sup-increments}
\end{equation}

\medskip
\noindent\textit{Step 1.}
Fix $m\in\N$.  It suffices first to work on $[0,m]$.  Choose
\begin{equation*}
 \ell>2(1+\beta),
 \qquad
 \delta_n:=2^{-\ell n},
 \qquad
 J_n:=\lceil m2^{\ell n}\rceil,
\end{equation*}
and set
\begin{equation*}
 x_{k,n}:=\frac{km}{J_n},\quad0\leq k\leq J_n,
 \qquad
 I_{j,n}:=(x_{j-1,n},x_{j,n}],\quad1\leq j\leq J_n.
\end{equation*}
Then $|I_{j,n}|\leq\delta_n$,
$J_n\leq C_m2^{\ell n}$, and $\log(J_n+1)\leq C_mn$.  Let
\begin{equation*}
 \xi_{j,n}
 :=D_n(x_{j,n})-D_n(x_{j-1,n})
 =\int_{(0,t]\times I_{j,n}\times(\varepsilon_{n+1},\varepsilon_n]}
 r\,\widetilde N(\dd s,\dd z,\dd r).
\end{equation*}
The $\xi_{j,n}$ are square-integrable and orthogonal because the deterministic spatial cells are disjoint.  The compensated-integral isometry and \eqref{eq:mean-measures} give
\begin{align*}
 \sum_{j=1}^{J_n}\Emu[\xi_{j,n}^2]
 &=c_\beta\Emu[Y_t((0,m])]
   \int_{\varepsilon_{n+1}}^{\varepsilon_n}r^{-\beta}\dd r
 \leq C_m2^{-n(1-\beta)}.
\end{align*}
Since $D_n(x_{k,n})=\sum_{j=1}^k\xi_{j,n}$, \eqref{eq:rademacher-menshov} yields
\begin{equation}
 \Emu\max_{0\leq k\leq J_n}|D_n(x_{k,n})|^2
 \leq C_mn^2 2^{-n(1-\beta)}.
 \label{eq:grid-maximal-bound}
\end{equation}
Choose $0<\kappa<(1-\beta)/2$.  By Markov's inequality and \eqref{eq:grid-maximal-bound},
\begin{equation*}
 \Pmu\!\left(
  \max_{0\leq k\leq J_n}|D_n(x_{k,n})|>2^{-\kappa n}
 \right)
 \leq C_mn^2 2^{-n(1-\beta-2\kappa)}.
\end{equation*}
The series on the right is summable.  Hence Borel--Cantelli and $\sum_n2^{-\kappa n}<\infty$ give
\begin{equation}
 \sum_{n=1}^\infty
 \max_{0\leq k\leq J_n}|D_n(x_{k,n})|<\infty
 \qquad\Pmu\text{-a.s.}
 \label{eq:grid-summability}
\end{equation}

\medskip
\noindent\textit{Step 2.}
Let $\Omega_L$ be the probability-one event in \eqref{eq:occupation-density}--\eqref{eq:occupation-consequences}, and set
\[
 \Omega_{m,M}
 :=\Omega_L\cap
 \left\{\sup_{z\in[-m,m]}L(t,z)\leq M\right\}.
\]
By \eqref{eq:occupation-consequences}, $\Pmu(\bigcup_{M\geq1}\Omega_{m,M})=1$.  For $1\leq j\leq J_n$, let
\begin{equation*}
 N_{j,n}(u)
 :=N((0,u]\times I_{j,n}\times(\varepsilon_{n+1},\varepsilon_n]),
 \qquad0\leq u\leq t.
\end{equation*}
Its continuous compensator is
\begin{equation*}
 \Lambda_{j,n}(u)
 =c_\beta\int_0^uX_s(I_{j,n})\dd s
   \int_{\varepsilon_{n+1}}^{\varepsilon_n}r^{-2-\beta}\dd r.
\end{equation*}
On $\Omega_{m,M}$,
\[
 Y_t(I_{j,n})\leq M\delta_n,
 \qquad
 \Lambda_{j,n}(t)\leq C_M2^{-n(\ell-1-\beta)}.
\]
Set $v_{n,M}:=C_M2^{-n(\ell-1-\beta)}$.  Pathwise,
\[
 \Omega_{m,M}\cap
 \left\{\max_{1\leq j\leq J_n}N_{j,n}(t)\geq2\right\}
 \subseteq
 \bigcup_{j=1}^{J_n}
 \{N_{j,n}(t)\geq2,\ \Lambda_{j,n}(t)\leq v_{n,M}\}.
\]
Thus the union bound and \eqref{eq:two-jump-bound} give
\begin{align*}
 &\Pmu\left(\Omega_{m,M}\cap
 \left\{\max_{1\leq j\leq J_n}N_{j,n}(t)\geq2\right\}\right)
 \leq CJ_nv_{n,M}^2
 \leq C_M2^{-n(\ell-2-2\beta)}.
\end{align*}
The series is summable because $\ell>2(1+\beta)$.  Borel--Cantelli, first on each $\Omega_{m,M}$ and then over $M$, shows that almost surely, eventually in $n$, every grid cell contains at most one raw jump with size in $(\varepsilon_{n+1},\varepsilon_n]$.

The raw-jump oscillation of $D_n$ within a cell is therefore eventually at most $\varepsilon_n=2^{-n}$.  On $\Omega_{m,M}$, the compensator oscillation over a cell is at most
\begin{equation*}
 c_\beta Y_t(I_{j,n})
 \int_{\varepsilon_{n+1}}^{\varepsilon_n}r^{-1-\beta}\dd r
 \leq C_M2^{-n(\ell-\beta)}.
\end{equation*}
Consequently, on $\Omega_{m,M}$, eventually in $n$,
\begin{equation*}
 \max_{1\leq j\leq J_n}\sup_{x,z\in I_{j,n}}|D_n(x)-D_n(z)|
 \leq 2^{-n}+C_M2^{-n(\ell-\beta)}.
\end{equation*}
The finitely many exceptional $n$ contribute a finite amount because each $D_n$ is bounded on compact intervals.  Hence
\begin{equation}
 \sum_{n=1}^\infty
 \max_{1\leq j\leq J_n}\sup_{x,z\in I_{j,n}}|D_n(x)-D_n(z)|<\infty
 \qquad\Pmu\text{-a.s.}
 \label{eq:cell-oscillation-summability}
\end{equation}

\medskip
\noindent\textit{Step 3.}
If $x\in I_{j,n}$, then
\[
 |D_n(x)|
 \leq|D_n(x_{j,n})|
 +\sup_{u,v\in I_{j,n}}|D_n(u)-D_n(v)|.
\]
Equations \eqref{eq:grid-summability} and \eqref{eq:cell-oscillation-summability} imply the positive-half version of \eqref{eq:summable-sup-increments}; the negative half is identical up to the sign in \eqref{eq:truncated-cumulative}.  Taking a countable intersection over $m$, define outside one null set
\[
 A_t(x):=\lim_{N\to\infty}A_t^{\varepsilon_N}(x),
 \qquad x\in\R,
\]
and set $A_t\equiv0$ otherwise.  Then, with probability one,
\begin{equation}
 A_t^{\varepsilon_N}\longrightarrow A_t
 \quad\text{uniformly on every compact interval.}
 \label{eq:compact-uniform-convergence}
\end{equation}
A locally uniform limit of c\`adl\`ag functions is c\`adl\`ag, and the construction also makes $A_t$ jointly measurable: for every $N$, the map $(\omega,x)\mapsto A_t^{\varepsilon_N}(\omega,x)$ is jointly measurable as the spatial cumulative function of two finite random measures, and $A_t$ is their pointwise limit on the measurable master event used above (and is set to zero on its complement).

\medskip
\noindent\textit{Step 4.}
Let $x\geq0$ be deterministic.  The martingale-measure isometry gives
\begin{align*}
 &\Emu\left|
 M_t(\1_{(0,x]})-M_t^{\varepsilon_n}(\1_{(0,x]})
 \right|^2\notag\\
 &\quad=c_\beta\Emu[Y_t((0,x])]
 \int_0^{\varepsilon_n}r^{-\beta}\dd r
 =\frac{c_\beta}{1-\beta}\Emu[Y_t((0,x])]
 \varepsilon_n^{1-\beta}\longrightarrow0.
\end{align*}
Thus $A_t^{\varepsilon_n}(x)\to M_t(\1_{(0,x]})$ in probability, whereas \eqref{eq:compact-uniform-convergence} gives almost-sure convergence to $A_t(x)$.  Uniqueness of limits in probability proves $A_t(x)=M_t(\1_{(0,x]})$ almost surely.  The same argument for $x<0$ gives $A_t(x)=-M_t(\1_{(x,0]})$ almost surely, establishing \eqref{eq:cumulative-modification}.

\medskip
\noindent\textit{Step 5.}
For $n\geq1$ and Borel $B\subset\R$, set
\[
 \nu_n(B)
 :=\int_{(0,t]\times B\times(\varepsilon_n,\infty)}
 r\,N(\dd s,\dd z,\dd r)
 -\int_{(0,t]\times B\times(\varepsilon_n,\infty)}
 r\,\widehat N(\dd s,\dd z,\dd r).
\]
By \eqref{eq:truncated-count-integrability} and \eqref{eq:truncated-compensator-integrability}, outside one null set every $\nu_n$ is a finite signed measure.  From \eqref{eq:truncated-cumulative},
\[
 A_t^{\varepsilon_n}(x)=
 \begin{cases}
  \nu_n((0,x]),&x\geq0,\\
  -\nu_n((x,0]),&x<0.
 \end{cases}
\]
The pathwise jump formula for the cumulative function of a finite signed measure gives, simultaneously for every $y\in\R$,
\begin{equation}
 A_t^{\varepsilon_n}(y)-A_t^{\varepsilon_n}(y-)=\nu_n(\{y\}).
 \label{eq:finite-signed-measure-jump}
\end{equation}
Moreover,
\begin{align*}
 \nu_n(\{y\})
 &=\int_{(0,t]\times\{y\}\times(\varepsilon_n,\infty)}
 r\,N(\dd s,\dd z,\dd r)\\
 &\quad-c_\beta Y_t(\{y\})
 \int_{\varepsilon_n}^\infty r^{-1-\beta}\dd r.
\end{align*}
Since $Y_t(\{y\})=0$ simultaneously for all $y$, with probability one,
\begin{equation}
 A_t^{\varepsilon_n}(y)-A_t^{\varepsilon_n}(y-)
 =\int_{(0,t]\times\{y\}\times(\varepsilon_n,\infty)}
 r\,N(\dd s,\dd z,\dd r),
 \quad n\geq1,\ y\in\R.
 \label{eq:truncated-spatial-jump}
\end{equation}
The simultaneity in $y$ follows from the pathwise identity
\eqref{eq:finite-signed-measure-jump}.

\medskip
\noindent\textit{Step 6.}
Intersect the probability-one events from Steps~3 and~5.  Fix $\omega$ in this intersection and $y\in\R$, and choose $m>|y|$.  Uniform convergence also controls the left limits:
\begin{equation*}
 |A_t^{\varepsilon_n}(y-)-A_t(y-)|
 \leq\sup_{x\in[-m,m]}|A_t^{\varepsilon_n}(x)-A_t(x)|\longrightarrow0.
\end{equation*}
Hence, by \eqref{eq:truncated-spatial-jump} and monotone convergence,
\begin{align}
 A_t(y)-A_t(y-)
 &=\lim_{n\to\infty}
 \int_{(0,t]\times\{y\}\times(\varepsilon_n,\infty)}
 r\,N(\dd s,\dd z,\dd r)\notag\\
 &=\int_{(0,t]\times\{y\}\times(0,\infty)}
 r\,N(\dd s,\dd z,\dd r).
 \label{eq:limit-spatial-jump}
\end{align}
The left-hand side is finite.  Since $y$ was arbitrary, the integral on the right is finite and \eqref{eq:limit-spatial-jump} holds simultaneously for every $y$, proving \eqref{eq:cumulative-fibre-jumps}.


\small
\begin{thebibliography}{99}

\bibitem{Daw93}
D.~A. Dawson,
\newblock Measure-valued Markov processes,
\newblock in \emph{Ecole d'Et\'e de Probabilit\'es de Saint-Flour XXI--1991}, Lecture Notes in Mathematics 1541, Springer, Berlin, 1993, pp.~1--260.

\bibitem{EKR91}
N.~El Karoui and S.~Roelly,
\newblock Martingale properties, explosion and L\'evy--Khintchine representation of a class of measure-valued branching processes,
\newblock \emph{Stochastic Process. Appl.} \textbf{38} (1991), 239--266.

\bibitem{Fle88}
K.~Fleischmann,
\newblock Critical behavior of some measure-valued processes,
\newblock \emph{Math. Nachr.} \textbf{135} (1988), 131--147.

\bibitem{FMW10}
K.~Fleischmann, L.~Mytnik and V.~Wachtel,
\newblock Optimal local H\"older index for density states of superprocesses with $(1+\beta)$-branching mechanism,
\newblock \emph{Ann. Probab.} \textbf{38} (2010), 1180--1220.

\bibitem{GS04}
I.~I. Gikhman and A.~V. Skorokhod,
\newblock \emph{The Theory of Stochastic Processes I},
\newblock Classics in Mathematics, Springer, Berlin, 2004.

\bibitem{LeGall16}
J.-F. Le Gall,
\newblock \emph{Brownian Motion, Martingales, and Stochastic Calculus},
\newblock Springer, 2016.

\bibitem{LM05}
J.-F. Le Gall and L.~Mytnik,
\newblock Stochastic integral representation and regularity of the density for the exit measure of super-Brownian motion,
\newblock \emph{Ann. Probab.} \textbf{33} (2005), 194--222.

\bibitem{Meaney07}
C.~Meaney,
\newblock Remarks on the Rademacher--Menshov theorem,
\newblock in A.~McIntosh and P.~Portal (eds.), \emph{Proceedings of the CMA/AMSI Research Symposium: Asymptotic Geometric Analysis, Harmonic Analysis and Related Topics}, Proceedings of the Centre for Mathematics and Its Applications, Australian National University, Vol.~42, 2007, pp.~100--110.

\bibitem{MP03}
L.~Mytnik and E.~Perkins,
\newblock Regularity and irregularity of $(1+\beta)$-stable super-Brownian motion,
\newblock \emph{Ann. Probab.} \textbf{31} (2003), 1413--1440.

\bibitem{MX04}
L.~Mytnik and K.~Xiang,
\newblock Tanaka formulae for $(\alpha,d,\beta)$-superprocesses,
\newblock \emph{J. Theoret. Probab.} \textbf{17} (2004), 483--502.

\bibitem{Sug89}
S.~Sugitani,
\newblock Some properties for the measure-valued branching diffusion processes,
\newblock \emph{J. Math. Soc. Japan} \textbf{41} (1989), 437--462.

\bibitem{Xiang05}
K.~Xiang,
\newblock On Tanaka formulae for $(\alpha,d,\beta)$-superprocesses,
\newblock \emph{Sci. China Ser. A Math.} \textbf{48} (2005), 1194--1208.

\end{thebibliography}
\end{document}